\documentclass[12pt]{amsart}
\usepackage{amsmath}
\usepackage{amsfonts,amssymb,amsthm}
\newcommand{\thm}[2]{\begin{#1} #2 \end{#1}}

\newtheorem{theorem}{Theorem}[section]

\newcommand{\integers}{{\bf Z}}

\begin{document}


\title{An inequality on Chebyshev polynomials}

\author{Igor Rivin}


\address{Department of Mathematics, Temple University, Philadelphia}

\curraddr{Mathematics Department, Princeton University}

\email{rivin@math.temple.edu}

\date{today}

\keywords{Chebyshev polynomials}

\subjclass{Primary 05C25, 05C20, 05C38, 60J10, 60F05, 42A05;
Secondary 22E27}

\begin{abstract}
We define a class of multivariate Laurent polynomials closely
related to Chebyshev polynomials, and prove the simple but
somewhat surprising (in view of the fact that the signs of the
coefficients of the Chebyshev polynomials themselves alternate)
result that their coefficients are non-negative.
\end{abstract}

\maketitle

\section*{Introduction}

Let $T_n(x)$ and $U_m(x)$ be the Chebyshev polynomials of the
first and second kinds, respectively. We define the following
Laurent polynomials:

\begin{eqnarray}
R_n(c; x_1, \dots, x_n) &=& T_n\left(\frac{c\sum_{i=1}^k \left[x_i
+ \frac{1}{x_i}\right]}{2 k}\right),\\
S_n(c; x_1, \dots, x_n) &=& U_n\left(\frac{c\sum_{i=1}^k \left[x_i
+ \frac{1}{x_i}\right]}{2 k}\right).
\end{eqnarray}
The function $R_n$ arises in the enumeration of conjugacy classes
in the free group on $k$ generators, more specifically, there is
the following result:
\begin{theorem}[\cite{rwalks}]
  \label{homoenum} The number of
cyclically reduced words of length $k$ in $F_r$ homologous to $e_1
[a_1] + \cdots + e_r [a_r]$ is equal to the coefficient of
$x_1^{e_1} \cdots x_r^{e_r}$ in
\begin{equation}
\label{genfn} 2\left(\sqrt{2r-1}\right)^k R_k({\frac{r}
{\sqrt{2r-1}}}; x_1, \dots, x_r)  + (r-1)[1 + (-1)^k]
\end{equation}
\end{theorem}

\medskip\noindent
{\em Remark.} The rescaled Chebyshev polynomial $T_k(a x)/a^k$ is
called the $k$-th Dickson polynomial $T_k(x, a)$ (see
\cite{schur}).

The coefficients of generating functions of combinatorial objects
are non-negative, so one is led to wonder for which values of $c$
are the coefficients of $R_n(c; x_1, \dots, x_k)$ and $S_n(c; x_1,
\dots, x_k)$ nonnegative. In this note we give an essentially
complete answer (see Theorems \ref{nonneg} and \ref{mnonneg}). We
also write down an explicit formula (Eq. (\ref{fullform})) for the
coefficients of $R_n$ and $S_n.$

\section{Some facts about Chebyshev polynomials}
\label{chebint}

The literature on Chebyshev polynomials is enormous; \cite{rivlin}
is a good to start. Here, we shall supply the barest essentials in
an effort to keep this paper self-contained.

There are a number of ways to define Chebyshev polynomials (almost
as many as there are of spelling their inventor's name). A
standard definition of the {\em Chebyshev polynomial of the first
kind} $T_n(x)$ is:

\begin{equation}
\label{def1} T_n(x) = \cos n \arccos x.
\end{equation}
In particular, $T_0(x) = 1,$ $T_1(x) = x.$ Using the identity
\begin{equation}
\label{cosid} \cos(x+y) + \cos(x-y) = 2 \cos x \cos y
\end{equation}
we immediately find the three-term recurrence for Chebyshev
polynomials:
\begin{equation}
\label{threerec} T_{n+1}(x) = 2 x T_n(x) - T_{n-1}(x).
\end{equation}
The definition of Eq. (\ref{def1}) can be used to give a ``closed
form'' used in Section \ref{homoesec}:
\begin{equation}
\label{sqrtdef} T_n(x) = {\frac12}\left[\left(x -
\sqrt{x^2-1}\right)^n + \left(x + \sqrt{x^2-1}\right)^n\right].
\end{equation}
Indeed, let $x = \cos \theta.$ Then
$$\left(x -
\sqrt{x^2-1}\right)^n = \exp(-i n \theta),$$ while
$$\left(x +
\sqrt{x^2-1}\right)^n = \exp(i n \theta),$$ so that
$${\frac12}\left(x - \sqrt{x^2-1}\right)^n + \left(x +
\sqrt{x^2-1}\right)^n = \Re \exp(i n \theta) = \cos n \theta.$$

We also define Chebyshev polynomials of the second kind $U_n(x)$,
which can again be defined in a number of ways, one of which is:
\begin{equation}
\label{derivdef} U_n(x) = {\frac1{n+1}}T_{n+1}^\prime(x).
\end{equation}
A simple manipulation shows that if we set $x = \cos \theta,$ as
before, then
\begin{equation}
\label{trigdef} U_n(x) = \frac{\sin (n+1) \theta}{\sin \theta}.
\end{equation}
In some ways, Schur's notation $\mathcal{ U}_n = U_{n-1}$ is
preferable. In any case, we have $U_0(x) = 1$, $U_1(x) = 2 x,$ and
otherwise the $U_n$ satisfy the same recurrence as the $T_n$, to
wit,
\begin{equation}
\label{trec2} U_{n+1}(x) = 2 x U_n(x) - U_{n-1}(x).
\end{equation}
From the recurrences, it is clear that for $f=T, U$, $f_n(-x) =
(-1)^n f(x),$ or, in other words, every second coefficient of
$T_n(x)$ and $U_n(x)$ vanishes. The remaining coefficients
alternate in sign; here is the explicit formula for the
coefficient $c_{n-2m}^{(n)}$ of $x^{n-2m}$ of $T_n(x):$
\begin{equation}
\label{coefform} c_{n-2m}^{(n)} = (-1)^m {\frac{n}{n-m}}
{\binom{n-m}{m}} 2^{n-2m-1}, \qquad m=0, 1, \ldots,
\left[{\frac{n}2}\right].
\end{equation}
This can be proved easily using Eq. (\ref{threerec}).

\section{Analysis of the functions $R_n$ and $S_n$.}
\label{genanal}

In view of the alternation of the coefficients, the appearance of
the Chebyshev polynomials as generating functions in Section
\ref{homoesec} seems a bit surprising, since combinatorial
generating functions have non-negative coefficients. Below we
state and prove a generalization. Remarkably, Theorems
\ref{nonneg} and \ref{mnonneg} do not seem to have been previously
noted.

\thm{theorem} { \label{nonneg} Let $c > 1.$ Then all the
coefficients of $R_n(c; x)$ are non-negative. Indeed the
coefficients of $x^n, x^{n-2}, \ldots, x^{-n + 2}, x^{-n}$ are
positive, while the other coefficients are zero. The same is true
of $S_n$ in place of $R_n.$ }

\begin{proof} Let $a_n^k$ be the coefficient of $x^k$ in
$U_n((c/2)(x+1/x)).$ The recurrence gives the following recurrence
for the $a_n^k:$
\begin{equation}
\label{newrec} a_{n+1}^k = c(a_n^{k-1}+a_n^{k+1}) - a_{n-1}^k.
\end{equation}
Now we shall show that the following always holds:

\begin{description}
\item[(a)] $a_n^k \geq 0$ (inequality being strict if and only if
$n-k$ is even).

\item[(b)] $a_n^k \geq \max(a_{n-1}^{k-1}, a_{n-1}^{k+1}),$ the
inequality strict, again, if and only if $n-k$ is even.

\item[(c)] $a_n^k \geq a_{n-2}^k$ (strictness as above).
\end{description}

The proof proceeds routinely by induction; first the induction
step (we assume throughout that $n-k$ is even; all the quantities
involved are obviously $0$ otherwise):

By induction $a_{n-1}^k < \min(a_n^{k-1}, a_n^{k+1}),$ so by the
recurrence \ref{newrec} it follows that $a_{n+1}^k >
\max(a_n^{k-1}, a_n^{k+1}).$ (a) and (c) follow immediately.

For the base case, we note that $a_0^0 = 1,$ while $a_1^1 =
a_1^{-1} = c > 1,$ and so the result for $U_n$ follows. Notice
that the above proof does {\em not} work for $T_n$, since the base
case fails. Indeed, if $b_n^k$ is the coefficient of $x^k$ in
$T_n((c/2)(x+1/x))$, then $b_0^0 = 1$, while $b_1^1 = c/2$, not
necessarily bigger than one. However, we can use the result for
$U_n$, together with the observation (which follows easily from
the addition formula for $\sin$) that
\begin{equation}
\label{moretrig} T_n(x) = {\frac{U_n(x) - U_{n-2}(x)}2}.
\end{equation}
Eq. (\ref{moretrig}) implies that $b_n^k = a_n^k - a_n^{k-2} > 0$,
by (c) above.
\end{proof}

The proof above goes through almost verbatim to show:

\thm{theorem} { \label{mnonneg} Let $c > 1.$ Then all the
coefficients of $R_n$ are non-negative. The same is true of $S_n$
in place of $R_n$ }

To complete the picture, we note that:

\thm{theorem} { \label{trivthm}
$$
R_n(1; x) = {\frac12}\left(x^n + {\frac1 {x^n}}\right).
$$
}

\begin{proof} Let $x = \exp i\theta.$ Then $1/2(x+1/x) = \cos \theta,$
and $R_n(1; x)=T_n(1/2(x+1/x)) = \cos n \theta = 1/2(x^n +
1/x^n).$
\end{proof}

\thm{remark} { For $c<-1$ it is true that all the coefficients of
$R_n(c; .)$ and $S_n(c; .)$ have the same sign, but the sign is
$(-1)^n.$ For $|c|<1,$ the result is completely false. For $c$
imaginary, the result is true. I am not sure what happens for
general complex $c$. }

By the formula (\ref{coefform}), we can write
\begin{equation}
\label{uniexp} T_n\left({\frac{c}
2}\left(x+{\frac{1}{x}}\right)\right) = {\frac12}
\sum_{m=0}^{\left[{\frac{n}{2}}\right]} (-1)^m {\frac{n}{n-m}}
{\binom{n-m}{m}} c^{n-2m} \left(x+{\frac{1}{x}}\right)^{n - 2m}.
\end{equation}
Noting that
\begin{equation}
\label{binom} \left(x+{\frac{1}{x}}\right)^k = \sum_{i=0}^k
{\binom{k} { i}} x^{k - 2i}
\end{equation}
we obtain the expansion
\begin{equation}
\label{fullform} R_n(c; x) = c^n \sum_{k=-n}^n x^k
\sum_{m=0}^{\left[{\frac{n}{2}}\right]}
\left(-{\frac{1}{c^2}}\right)^m {\frac{n}{n-m}} {\binom{n-m}{m}}
{\binom{n-2m}{(n-2m-k)/2}},
\end{equation}
where it is understood that $\binom{a} { b}$ is $0$ if $b<0,$ or
$b > a$, or $b \notin \integers.$ A similar formula for $S_n$ can
be written down by using Eq. (\ref{moretrig}).

\bibliographystyle{amsplain}

\end{document}